\pdfoutput=1
\RequirePackage{ifpdf}
\ifpdf 
\documentclass[pdftex]{sigma}
\else
\documentclass{sigma}
\fi

\def\D{\mathcal D}
\def\F{\mathcal F}
\def\K{\mathcal K}

\def\L{\mathcal L}
\def\H{\mathcal H}
\def\O{\mathcal O}
\def\T{\mathcal T}

\def\1{\mathbf 1}
\def\M{\mathcal M}

\def\QQ{\mathbb Q}
\def\ZZ{\mathbb Z}
\def\CC{\mathbb C}

\def\Res{\operatorname{Res}}
\def\hat{\widehat}
\def\tilde{\widetilde}
\def\p{\partial}
\def\a{\alpha}

\def\f{{\mathbf f}}
\def\g{{\mathbf g}}
\def\t{{\mathbf t}}

\def\pp{{\mathbf p}}

\def\ll{{\mathbf l}}
\def\vv{{\mathbf v}}
\def\uu{{\mathbf u}}
\def\ww{{\mathbf w}}

\def\gs{\sigma}

\def\h{\hbar}

\def\lan{\langle}
\def\ran{\rangle}
\def\str{\operatorname{str}}

\def\ev{\operatorname{ev}}

\def\td{\operatorname{td}}
\def\ch{\operatorname{ch}}

\def\eu{\operatorname{eu}}

\def\Eu{\operatorname{Eu}}
\def\tr{\operatorname{tr}}
\def\str{\operatorname{str}}

\def\square{\Box}
\def\und{\underline}
\renewcommand{\Delta}{\triangle}

\begin{document}
\allowdisplaybreaks

\newcommand{\arXivNumber}{1710.02376}

\renewcommand{\thefootnote}{}

\renewcommand{\PaperNumber}{031}

\FirstPageHeading

\ShortArticleName{Quantum Hirzebruch--Riemann--Roch in Genus 0}

\ArticleName{Permutation-Equivariant Quantum K-Theory~X. \\ Quantum Hirzebruch--Riemann--Roch in Genus 0\footnote{This paper is a~contribution to the Special Issue on Algebra, Topology, and Dynamics in Interaction in honor of Dmitry Fuchs. The full collection is available at \href{https://www.emis.de/journals/SIGMA/Fuchs.html}{https://www.emis.de/journals/SIGMA/Fuchs.html}}}

\Author{Alexander GIVENTAL}
\AuthorNameForHeading{Alexander Givental}
\Address{Department of Mathemetics, UC Berkeley, CA 94720, USA}
\Email{\href{mailto:givental@math.berkeley.edu}{givental@math.berkeley.edu}}
\URLaddress{\url{https://math.berkeley.edu/~giventh/}}

\ArticleDates{Received September 28, 2019, in final form April 13, 2020; Published online April 22, 2020}

\Abstract{We extract genus $0$ consequences of the all genera quantum HRR formula proved in Part~IX. This includes re-proving and generalizing the adelic characterization of genus~$0$ quantum K-theory found in [Givental A., Tonita V., in Symplectic, {P}oisson, and Noncommutative Geometry, \textit{Math. Sci. Res. Inst. Publ.}, Vol.~62, Cambridge University Press, New York, 2014, 43--91]. Extending some results of Part VIII, we derive the invariance of a certain variety (the ``big J-function''), constructed from the genus~$0$ descendant potential of permutation-equivariant quantum K-theory, under the action of certain finite difference operators in Novikov's variables, apply this to reconstructing the whole variety from one point on it, and give an explicit description of it in the case of the point target space.}


\Classification{14N35}

\renewcommand{\thefootnote}{\arabic{footnote}}
\setcounter{footnote}{0}

\rightline{\em To my K-theory teacher on his 80th B-day}

\section{Adelic characterization}

Let $\F_X$ denote the genus-0 descendant potential of permutation-equivariant
quantum K-theory on a compact K\"ahler manifold $X$, or more precisely, its {\em dilaton-shifted} version. By definition (see Part IX)
\[ \F_X(\vv+\t) := \sum_{\ll, d} \frac{Q^d}{\prod_r l_r!} \lan \dots, \t_1, \dots; \dots ; \dots, \t_k, \dots; \dots \ran_{0,\ll,d}.\]
Here $\t=(\t_1,\dots, \t_k,\dots)$ is a sequence of Laurent polynomials in $q$ with vector coefficients in $K:=K^0(X)\otimes \Lambda$, $\vv :=(1-q) (\1,\1,\dots)$ is the {\em dilaton vector}, and $\1$ is the unit in $K^0(X)$. The {\em ground ring} $\Lambda$, must contain Novikov's variables $Q$, and be a $\lambda$-algebra, i.e., equipped with Adams' operations $\Psi^m$, $m=1,2,\dots$, $\Psi^1=\operatorname{id}$, which act on Novikov's variables by $\Psi^m\big(Q^d\big)=Q^{md}$. We assume that $\Psi^m$ with $m>1$ increase the descending filtration of $\Lambda$ by the powers $\Lambda_{+}^d$ of a certain ideal $\Lambda_{+}$ containing all Novikov's monomials $Q^d$ with $d\neq 0$. The correlators $\lan \dots \ran_{g,\ll,d}$ are defined in terms of K-theory on the moduli space $X_{g,n,d}$ of genus-$g$ degree-$d$ stable maps to $X$ of compact complex nodal curves with $n$ marked points. Let $h$ denote the automorphism of $X_{g,n,d}$ induced by the renumbering of the $n$ marked points with the cycle structure
determined by a partition $\ll=(l_1,\dots,l_r,\dots)$, where $l_r$ stands for the number of cycles of length $r=1,2,\dots$. The correlator
\[ \lan \uu_1,\dots, \uu_{l_1}; \vv_1, \dots, \vv_{l_2}; \dots ; \ww_1, \dots, \ww_{l_r}; \dots \ran_{g,\ll,d},\]
where $\uu_i$, $\vv_j$, $\ww_k$ are $K$-valued Laurent polynomials of $q$, is defined as follows. To such a~Laurent polynomial, say, $\ww(q)=\sum_m \phi_m q^m$, $\phi_m\in K^0(X)\otimes \Lambda$, we associate an $h$-equivariant vector bundle~$W$ on $X_{g,n,d}$:
\[ W = \bigotimes_{\a=1}^r \sum_m \ev_{\sigma_\a}^*(\phi_m) L_{\sigma_\a}^{\otimes m},\]
where $\sigma_1,\dots,\sigma_r$ are indices of marked points cyclically permuted by $h$, $\ev_{\sigma}\colon X_{g,n,d}\to X$ is the evaluation map at the marked point $\sigma$, and $L_{\sigma}$ is the line orbibundle over $X_{g,n,d}$ formed by the cotangent lines to the curves at the $\sigma$th marked point.
With this notation, the above correlator equals
\[ \left(\prod_{r=1}^{\infty}\frac{1}{ r^{l_r}}\right) \ \str_h H^*\left(X_{g,n,d}; \O_{g,n,d}\bigotimes_{i=1}^{l_1} U_i \bigotimes_{j=1}^{l_2} V_j \cdots \bigotimes_{k=1}^{l_r}W_k \cdots \right) ,\]
where $\O_{g,n,d}$ is the virtual structure sheaf introduced by Y.-P. Lee \cite{YPLee}, and
\[\str_h H^* := \sum (-1)^k\tr_h H^k\] denotes the alternated sum of traces of the operator $h$ on the graded space.

With the function $\F_X$ we associate the ``big J-function'' which to a
sequence $\t = (\t_1,\t_2,\dots)$ of $K$-valued Laurent polynomials in $q$
associates the following sequence $\f = (\f_1,\f_2,\dots)$ of $K$-valued rational functions\footnote{We remind that by a rational function, Laurent polynomial, etc. in $q$ we mean a $Q$-series whose coefficients are such rational functions, Laurent polynomials, etc.} of $q$:
\begin{gather*}
 \f_1 = (1-q)\1 + \t_1(q) + \sum_{\ll, d, \a} \frac{Q^d\phi_\a}{\prod_s {l_s!}}
 \left\langle \frac{\phi^{\a}}{1-qL}, \t_1,\dots; \t_2, \dots \right\rangle_{d, \ll+\1_1}, \\
 \cdots\cdots\cdots\cdots\cdots\cdots\cdots\cdots\cdots\cdots\cdots\cdots\cdots\cdots\cdots\cdots\cdots\cdots \cdots\cdots\cdots\\
 \f_r = (1-q)\1 + \t_r(q) + \sum_{\ll, d, \a} \frac{Q^d\phi_\a}{\prod_s {l_s!}}\left\langle \frac{\phi^{\a}}{1-qL}, \t_k, \dots ; \t_{2k}, \dots \right\rangle_{d, \ll+\1_1}, \\
 \cdots\cdots\cdots\cdots\cdots\cdots\cdots\cdots\cdots\cdots\cdots\cdots\cdots\cdots\cdots\cdots\cdots\cdots \cdots\cdots\cdots \end{gather*}
Here the cycle structure $\ll+\1_1$ is obtained from $\ll$ by adding one cycle of length $1$, and $L$ is the line (orbi)bundle on the moduli
spaces formed by the cotangent lines to the curve at the 1st (unpermutted) marked point. Note that on each particular moduli space, due to its orbifold nature, $L$ is not unipotent, but nevertheless satisfies $P(L)=0$, where
$P$ is some Laurent polynomial (with zeroes -- in the non-equivariant setting -- at some roots of unity). Since $P\big(q^{-1}\big)-P(L)=(1-qL) \Phi(q,L)$, where~$\Phi$ is some Laurent polynomial in two variables, we have
$1/(1-qL)=\Phi(q,L)/P(q)$. This shows that each correlator in the series is a rational function of $q$ (with poles at some roots of unity), and -- obviously
from the definition -- having no pole at $q=0$ and vanishing at $q=\infty$.
The part $(1-q)\1+\t_r(q)$ is, in the contrary, a Laurent polynomial in~$q$,
i.e., can have no poles except $q=0, \infty$.

The range ($=$ image) of the J-function, viewed as a map $\t \mapsto \f (\t)$, is a (formal, infinite dimensional) subvariety in the space $\K^{\infty}$ of sequences $\f = (\f_1,\f_2,\dots)$ of $K$-valued rational functions in~$q$ allowed to have poles at $q=0,\infty$, or at roots of unity. We denote this range by~$\L_X$. It can also be viewed as the graph of (the formal germ at
$(1-q)(\1,\1,\dots)$ of) the map
$\K^{\infty}_{+}\ni \t \mapsto \f - (1-q)(\1,\1,\dots) \in \K^{\infty}_{-}$,
where $\K^{\infty}_{+}$ consists of sequences of $K$-valued Laurent polynomials,
and $\K^{\infty}_{-}$ consists of sequences of rational functions
having no pole at $q=0$ and vanishing at $q=\infty$.

Yet another way to describe $\L_X$ is to consider $\F_X$ as a family of functions of the first input~$\t_1$ depending on $\t_2, \t_3,\dots$ as parameters. Then the map $\t_1 \mapsto \f_1(\t_1, \t_2,\t_3,\dots)$ para\-meterizes the graph of the differential ${\rm d}_{t_1}\F_X$ for given values of the parameters $\t_2,\t_3,\dots$, and $\f_r:=\f_1(\t_r,\t_{2r},\t_{3r},\dots)$.
The sequences $\f=(\f_1,\f_2,\dots )$ of the differentials lie in the cotangent bundle~$T^*\K_{+}^{\infty}$ of the subspace $\K_{+}^{\infty}\subset \K^{\infty}$. The cotangent bundle space is identified with the whole of~$\K^{\infty}$ by means of the polarization $\K^{\infty}=\K^{\infty}_{+}\oplus\K^{\infty}_{-}$, which is Lagrangian with respect to the $\Lambda$-valued symplectic form
\begin{gather*} \Omega^{\infty}(\f,\g) = \sum_{r=1}^{\infty} \frac{\Psi^r}{r} \Omega (\f_r,\g_r),\ \\
 \Omega(f,g) = - [\Res_{q=0}+\Res_{q=\infty}] \big(f\big(q^{-1}\big), g(q)\big) \frac{{\rm d}q}{q}.\end{gather*}
Here $(a, b):=\chi (X; a\otimes b) = \int_X \td(T_X)\ch(a)\ch(b)$ is the
K-theoretic Poincar\'e pairing on $K^0(X)$.
By definition, the negative space $\K_{-}^{\infty}$ of the polarization consists of sequences $\f = (\f_1,\f_2,\dots)$ satisfying $\f_r(\infty)=0$, $\f_r(0)\neq \infty$.

\begin{remark} Note that the graph of each $\f_r$ for fixed values of the parameters $\t_{2r}, \t_{3r}, \dots$ is a~Lagrangian variety in $(\K, \Omega)$ (where $\K$ is the space of $K$-valued Laurent polynomials in~$q$), and the whole family of such Lagrangian varieties for, say, $r=1$ determines $\L_X$. However, a value of~$\f_1$, though ``knows'' its input $\t_1$, does not ``remember'' the values of the parameters $\t_2,\t_3,\dots$. One reason to introduce $\L_X$ parameterized by all the $\f_r$ is that from a point $\f=(\f_1,\f_2,\dots)\in \L_X$, the corresponding argument $\t=(\t_1,\t_2,\dots) \in \K^{\infty}_{+}$ is reconstructed by projection along $\K^{\infty}_{-}$. On the other hand, there is no reason for $\L_X$ to be Lagrangian in $(\K^{\infty},\Omega^{\infty})$.
\end{remark}

We will give an {\em adelic characterization} of $\L_X$ in terms of {\em fake}
genus-0 K-theoretic GW-invariants of $X$. The fake holomorphic Euler characteristic of a bundle $V$ over $X_{g,n,d}$ is defined by the right hand side of the Hirzebruch-RR formula:
\[ \chi^{\rm fake}(X_{g,n,d}; V):= \int_{[X_{g,n,d}]} \td (T_{X_{g,n,d}}) \ch (V),\]
where $[X_{g,n,d}]$ is the virtual fundamental class, and $T_{X_{g,n,d}}$ is the virtual tangent bundle. The genus-0 descendant potential $\F_X^{\rm fake}$ is defined as a function of $t(q)=\sum_{m\geq 0} \phi_m (q-1)^m$ by
\[ \F_X^{\rm fake}(t)=\sum_{n,d}\frac{Q^d}{n!} \lan t(L),\dots,t(L)\ran_{0,n,d}^{\rm fake},\]
where $\big\lan \phi_{\a_1} L^{m_1},\dots, \phi_{\a_n} L^{m_n}\big\ran_{0,n,d}:=\chi^{\rm fake}\big(X_{0,n,d};\otimes_{i=1}^n\ev_i^*(\phi_{\a_i}) L_i^{m_i}\big)$.

The graph of the differential ${\rm d}\F_X^{\rm fake}$ is identified
with a Lagrangian variety $\L_X^{\rm fake}\!\!\subset\!\! \big(\K^{\rm fake}\!,\Omega^{\rm fake}\big)$. By $\K^{\rm fake}$ we denote the space of $K$-valued Laurent series in $q-1$. It is equipped with the symplectic form
\[ \Omega^{\rm fake}(f,g):=\Res_{q=1} \big(f\big(q^{-1}\big),g(q)\big) \frac{{\rm d}q}{q},\]
Lagrangian polarization $\K^{\rm fake}_{\pm}$, where $\K^{\rm fake}_{+}$ consists of power series, and $\K^{\rm fake}_{-}$ consists of the principal parts of the Laurent series, as well as with the {\em dilaton vector} $(1-q)\1$.
Explicitly, $\L_X^{\rm fake}$ consists of vector-valued Laurent series in $q-1$ of the form
\[ 1-q + t(q) + \sum_{\a,n,d} \frac{Q^d\phi_{\a}}{n!}\left\lan \frac{\phi^{\a}}{1-qL}, t(L),\dots, t(L)\right\ran^{\rm fake}_{0,n+1,d}.\]
In fact $\L_X^{\rm fake}$ is an {\em overruled Lagrangian cone}, i.e., each tangent space $T$ of $\L_X^{\rm fake}$ contain $(1-q)T \subset \L_X^{\rm fake}$ and is tangent to $\L_X^{\rm fake}$ everywhere along $(1-q)T$. The cone $\L_X^{\rm fake}$ can be explicitly described in terms of cohomological GW-theory of $X$ (see~\cite{Co, GiF, GiTo} or Part~IX).

\begin{theorem}\label{Theorem1} Let $\f = (\f_1,\f_2,\dots, \f_r,\dots)$ be a sequence of rational functions of $q$ with values in $K=K^0(X)\otimes \Lambda$ considered as an element of the symplectic loop space
 $(\K^{\infty}, \Omega^{\infty})$. This sequence represents a point in the range
 $\L_X \subset \K^{\infty}$ of the ``big J-function'' of permutation-equivariant quantum K-theory of $X$ if and only if it satisfies the following
 criteria:
\begin{enumerate}\itemsep=0pt
\item[$(i)$] For each $r=1,2,\dots$, the Laurent series expansion $\f_r^{(1)}$ of $\f_r(q)$ near $q=1$ lies in the Lagrangian submanifold $\L_X^{\rm fake} \subset \big(\K^{\rm fake}, \Omega^{\rm fake}\big)$ representing the graph of the differential of the genus $0$ descendant potential $\F_X^{\rm fake}$ of fake quantum K-theory of~$X$.
\item[$(ii)$] For each $r$, and each root of unity $\zeta\neq 1$, the Laurent series expansion $\f_r^{(\zeta)}$ of $\f_r\big(q^{1/m}/\zeta\big)$ near $q=1$, where~$m$ is the primitive order of~$\zeta$, lies in the Lagrangian subspace
 \[ \Delta_{\zeta} \Psi^m\big(T_{\f_{rm}^{(1)}}\L_X^{\rm fake}\big) \otimes_{\Psi^m(\Lambda)}\Lambda,\]
 where the operator of multiplication $\Delta_{\zeta}$ is given by the formula
 \[ \Delta_{\zeta}:= \exp\left\{\sum_{k>0}\left(\frac{\Psi^{k}(T^*_X-1)}{k(1-\zeta^{-k}q^{k/m})} -\frac{\Psi^{km}(T^*_X-1)}{k(1-q^{km})}\right)\right\}.\] \end{enumerate}
\end{theorem}

Theorem~\ref{Theorem1} generalizes several previously known results. Setting $\t_k=0$ for all $k>1$ we obtain the result of \cite{GiTo} which characterizes the version of genus-0 quantum K-theory which does not involve permutations of the marked points. Setting all $\t_k$ with $k>1$ equal to each other, we obtain an adelic characterization of permutation-equivariant J-function studied in Part VII and
Part VIII \cite{GiVII,GiVIII}.

In fact one can prove Theorem~\ref{Theorem1} by generalizing the arguments in the genus-$0$ paper \cite{GiTo} by the author and Valentin Tonita. Our goal here is, however,
is to derive this genus-$0$ result from the all-genera quantum HRR theorem established in Part IX of this series. That all-genera result is a formula expressing a quantum state (in the Fock space of $(\K^{\infty},\Omega^{\infty})$), which represents the appropriate total descendant potential of the theory, in terms of its adelic counterpart described via the ``fake'' quantum K-theory. What seems peculiar at the first sight is that Theorem~\ref{Theorem1} above is not an {\em adelic formula}, but merely an {\em adelic characterization} of points on $\L_X$ (in terms of their localizations at the roots of unity). How can the quasi-classical limit at $\h=0$ (and this is how the genus-$0$ information is supposed to be extracted from the quantum-mechanical description of the all-genera data) fail to yield a formula too? The answer to this puzzle lies in the unusual feature of the quantum-mechanical formalism of the permutation-equivariant quantum K-theory, where the Adams operations act non-trivially on the Planck constant: $\Psi^r\big(\h^{\pm 1}\big)=\h^{\pm r}$. Consequently, the quasi-classical limit loses its familiar form of a single Lagrangian variety, and acquires a more intriguing structure.

In Sections~\ref{section2}--\ref{section4}, we give several aforementioned applications of Theorem~\ref{Theorem1}, and postpone its derivation from the all-genera quantum HRR theorem till the last
section.

\section[Example: $X=pt$]{Example: $\boldsymbol{X=pt}$}\label{section2}

In this case, $\K^{\infty}$ consists of sequences $\f=(\f_1,\dots,\f_r,\dots)$
of scalar (i.e., $\Lambda$-valued) rational functions of $q$ with poles at $q=0,\infty$ or roots of unity, and $\K_{+}^{\infty}$ consists of sequences
$\t = (\t_1,\dots,\t_r,\dots )$ of scalar Laurent polynomials in $q$.

\begin{theorem}\label{Theorem2} The range $\L_{pt}\subset \K^{\infty}$ of the ``big J-function'' consists of sequences $\f =(\f_1,\dots,\f_r,\allowbreak \dots) \in \K^{\infty}$ of the form
 \begin{gather*} \f_1(q) = (1-q) \exp\left\{ \sum_{k>0} \Psi^k(\tau_k)/k\big(1-q^k\big)\right\} \t_1\left(q,q^{-1}\right), \\
 \cdots\cdots\cdots\cdots\cdots\cdots\cdots\cdots\cdots\cdots\cdots\cdots\cdots\cdots\cdots\cdots\cdots \\
 \f_r(q) =(1-q) \exp\left\{ \sum_{k>0} \Psi^k(\tau_{kr})/k\big(1-q^k\big)\right\} \t_r\left(q,q^{-1}\right), \\
 \cdots\cdots\cdots\cdots\cdots\cdots\cdots\cdots\cdots\cdots\cdots\cdots\cdots\cdots\cdots\cdots\cdots \end{gather*}
 where $\t_k \in \K_{+}$ are scalar Laurent polynomials $\Lambda_{+}$-close to $1$, and $\tau_k\in \Lambda_{+}$ are arbitrary scalar parameters.
 \end{theorem}

\begin{proof}We derive this by verifying criteria (i), (ii) of Theorem~\ref{Theorem1}.

It is known (see, for instance, \cite{CGL}) that the $\L^{\rm fake}_{pt}$ consists of Laurent series of the form
\[ (1-q) {\rm e}^{\tau/(1-q)}\ t(q-1),\]
where $t$ is a power series in $q-1$ which is $\Lambda_{+}$-close to $1$, and
$\tau\in \Lambda_{+}$ is an arbitrary scalar.

First, we expand $\f_r$ into Laurent series $\f_r^{(1)}$ in $q-1$:
\[ \f_r^{(1)}= (1-q){\rm e}^{\T_r/(1-q)}\ t_r(q-1), \qquad \text{where}\qquad \T_r=\sum_{k>0}\Psi^k(\tau_{kr})/k^2,\]
and $t_r$ are some power series in $q-1$. We use here that $1/\big(1-q^k\big)=1/k(1-q)+\O(1)$. Thus $\f_r^{(1)}\in \L_{pt}^{\rm fake}$, and
criterion (i) is fulfilled.

Note that the tangent space to $\L_{pt}^{\rm fake}$ at the point $\f_r^{(1)}$ has the form ${\rm e}^{\T_r/(1-q)} \K_{+}^{\rm fake}$.

In the case $X=pt$, the operator
$\Delta_{\zeta}$ consists in multiplication by an invertible power series in $q-1$. Therefore, for a primitive $m$th root of unity $\zeta$,
\[ \Delta_{\zeta} \Psi^m\big(T_{\f_{rm}^{(1)}}\L_{pt}^{\rm fake}\big)\otimes_{\Psi^m(\Lambda)}\Lambda = {\rm e}^{\Psi^m[\T_{rm}/(1-q)]}\K^{\rm fake}_{+}.\]
We want to check that this subspace contains the Laurent series expansion $\f_r^{(\zeta)}$ of $\f_r\big(q^{1/m}/\zeta\big)$ near $q=1$. We have
\[ \f_r^{(\zeta)}=\big(1-q^{1/m}/\zeta\big)
\exp\left\{ \sum_{k>0} \Psi^k(\tau_{kr})/k\big(1-q^{k/m}/\zeta^k\big)\right\} \t_r\big(q^{1/m}/\zeta, q^{-1/m}\zeta\big).\]
The terms in the exponent which have pole at $q=1$ come from the values of $k$
divisible by $m$. For $k=lm$ we have
\[ \sum_{l>0} \frac{\Psi^{lm}(\tau_{lmr})}{lm(1-q^l)} = \Psi^m\left(\frac{\T_{rm}}{1-q}\right) + \text{terms regular at $q=1$}.\]
Thus, criterion (ii) is also fulfilled.

Thus, the family of sequences $\f=(\f_1,\dots,\f_r,\dots)$ described in Theorem~\ref{Theorem2} lie in $\L_{pt}$. Note that at $\tau = 0$, this family contains the $\Lambda_{+}$-neighborhood of the dilaton vector in $\K_{+}^{\infty}$. At $\t=\1$, the projection of $\f$ to $\K^{\infty}_{+}$ along $\K^{\infty}_{-}$ is
$[\f]_{+}=(1-q+\tau_1, \dots, 1-q+\tau_r,\dots)$. Thus, the tangent space at the dilaton point to the entire domain of $\F_{pt}$ is covered by our family. It follows now from the formal implicit function theorem, that the family parameterizes the whole of $\L_{pt}$.
\end{proof}

\section[$\D_q$-symmetries]{$\boldsymbol{\D_q}$-symmetries}\label{section3}

We generalize here results of \cite{GiE, GiTo}, and Part VIII about symmetries of genus-0 quantum K-theory induced by finite-difference operators in Novikov's variables.

Pick an integer basis $p_1,\dots, p_s$ in $H_2(X;\QQ)$ with respect to which all degrees of holomorphic curves in $X$ are expressed by vectors $d=(d_1,\dots,d_s)$ with non-negative components~$d_i$. Let~$P_i$ be the line bundle over~$X$ whose 1st Chern class is~$-p_i$. Let~$\D_q$ denote the algebra of finite difference operators in Novikov's variables. By definition it consists of non-commutative polynomial expressions formed from multiplication operators~$Q_i$, translation operators $q^{Q_i\p_{Q_i}}$, and can have Laurent polynomials in~$q$ in the role of coefficients. We make $\D_q$ act on the space $\K$ of $K$-valued rational functions of~$q$ so that~$Q_i$ act naturally as multiplication by $Q_i\in \Lambda_{+}$, while the translation operators act on
$K$-valued functions of $Q$ by $P_iq^{Q_i\p_{Q_i}}$, i.e., by
\[ \F (Q_1,\dots, Q_s) \mapsto P_i \F (\dots, Q_{i-1}, qQ_i, Q_{i+1},\dots
).\]
For convergence purposes we will further assume that our operators $D\big(Pq^{Q\p_Q},Q,q\big)$ have ``small free terms'', i.e., $D(1,0,q)\in \Lambda_{+}\big[q,q^{-1}\big]$.

\begin{theorem}\label{Theorem3} Let $D=(D_1,\dots,D_k,\dots)$ be a sequence of finite
 difference operators $D_k\big(Pq^{Q\p_Q},\linebreak Q, q\big)$. Then the following transformation
 on the space $\K^{\infty}$ of sequences $\f=(\f_1,\dots,\f_r,\dots)$ of vector-valued rational function preserves $\L_X$:
\[ \f_r \mapsto \exp\left\{ \sum_{k>0}
 \Psi^k \big(D_{kr}\big(Pq^{kQ\p_Q},Q,q\big)\big)/k\big(1-q^k\big)\right\} \f_r,\qquad r=1,2,\dots .\]
\end{theorem}

\begin{remark} Recall that Adams' operations $\Psi^k$ act naturally on
$K^0(X)$, act of functions of $q$ by $\Psi^k(q)=q^k$, act through the $\lambda$-algebra structure on the coefficient ring $\Lambda$, and in particular by $\Psi^k\big(Q^d\big)=Q^{kd}$ on Novikov's variables. Consequently, multiplication by $P_i$ hidden in $\Psi^k \big(D\big(Pq^{Q\p_Q},Q,q\big)\big)$ turns after the application of $\Psi^k$ into multiplication by $P_i^k$. However,
\[ \Psi^k \big(q^{Q_i\p_{Q_i}} Q^d\big) =\Psi^k\big(q^{d_i} Q^d\big) = q^{kd_i} Q^{kd} = q^{Q_i\p_{Q_i}} \Psi^k\big(Q^d\big),\]
i.e., $\Psi^k \big(q^Q\p_Q\big) =q^{Q\p_Q}$. Therefore $\Psi^k\big(Pq^{kQ\p_Q}\big)=\big(Pq^{Q\p_Q}\big)^k$. Thus, the operators in the exponent act as legitimate finite difference operators in our representation, i.e., as combinations of
multiplications by~$Q_i$ and twisted translation operators $P_iq^{Q_i\p_{Q_i}}$.
\end{remark}

\begin{proof} Suppose $\f \in \L_X$. Then: (i)~expansion $\f_r^{(1)}$ of $\f_r(q)$ near $q=1$ lie in $\L_X^{\rm fake}$, and
(ii)~expansion $\f_r^{(\zeta)}$ of $\f_r\big(q^{1/m}/\zeta\big)$ near $q=1$, where $\zeta\neq 1$ is a primitive $m$th root of unity, lie in $\Delta_{\zeta} \Psi^m\big(T_{\f_{rm}^{(1)}}\L_X^{\rm fake}\big) \otimes_{\Psi^m(\Lambda)}\Lambda$.

It is known (see, for instance, \cite{Co, GiTo}, or Section~3 in Part~IX), that
$\L_X^{\rm fake}$ is an overruled Lagrangian cone in the symplectic loop space~$\K^{\rm fake}$, and that it is obtained by a certain loop group transformation $\Delta\colon \H\to \K^{\rm fake}$ from the overruled Lagrangian cone $\L_X^H$ of quantum cohomology theory on $X$. Recall that $\L_X^H$ lies in the appropriate symplectic loop space $\H=H^{\rm even}(X; \Lambda ((z)))$, where $z=\log q$, and represents the graph of the differential of the genus~0 descendant potential of cohomological GW-theory of~$X$. It follows (see~\cite{GiE}) from the divisor equations of quantum cohomology theory, that $\L_X^H$ is invariant under the action of operators $f \mapsto {\rm e}^{D/z}f$ where $D(zQ\p_Q-p, Q,z)$ is any (pseudo) differential operator (with ``small free term'') in Novikov's variables. Consequently, since the loop group transformation~$\Delta$ does not depend on~$Q$, the cone $\L_X^{\rm fake}$ is similarly invariant under the action of operators
$f\mapsto {\rm e}^{D/(1-q)}f$, where $D((\log q)Q\p_Q-p, Q, q-1)$ is any differential operator.\footnote{Note that $\log q = \log (1-(1-q))=-\sum\limits_{k>0} (1-q)^k/k$ is a legitimate coefficient in $\K^{\rm fake}$ (though not in $\K$), and $p_i=-\log P_i$ is well-defined under the Chern character identification of $K^0(X)\otimes \QQ$ with $H^{\rm even}(X;\QQ)$.} The operator in the exponent here is required to have at most the 1st order pole at $q=1$. Besides, the tangent spaces $T_f\L_X^{\rm fake}$ and the {\em ruling spaces} $(1-q)T_f\L_X^{\rm fake} \subset \L_X^{\rm fake}$ are $\D$-modules, i.e., stay invariant under the operators $D$, and hence under the operators ${\rm e}^D$ (with no pole at $q=1$). In particular, since $Pq^{Q\p_Q} ={\rm e}^{(\log q) Q\p_Q-p}$, the invariance properties hold true for finite difference operators.

\looseness=-1 Now let $\g =(\g_1,\dots,\g_r,\dots)$ be the result of the transformation described in the theorem and applied to $\f$. We intend to check that $\g$ satisfies the criteria (i), (ii) of Theorem~\ref{Theorem1}. Since $1/\big(1-q^k\big)$ has the 1st order pole at $q=1$, it follows that $\g_r^{(1)}\in \L_X^{\rm fake}$. Furthermore, according to the invariance properties described above, the tangent space to $\L_X^{\rm fake}$ at $\f_{rm}^{(1)}$ is transformed to another tangent space (at $\g_{rm}^{(1)}$) determined by the polar part of the operator in the exponent. Therefore
\[ T_{\g_{rm}^{(1)}}\L_X^{\rm fake} = \exp\left\{\sum_{k>0} \Psi^k\big(D_{krm}\big(Pq^{kQ\p_Q},Q,q\big)\big)/k^2(1-q)\right\} T_{\f_{rm}^{(1)}}\L_X^{\rm fake} .\]
On the other hand,
\begin{gather*} \g_r^{(\zeta)} = \exp\left\{\sum_{k>0} \Psi^k \big(D_{kr}\big(P\big(q^{1/m}/\zeta\big)^{kQ\p_Q},Q,q^{1/m}/\zeta\big)\big)/k\big(1-q^{k/m}/\zeta^k\big)\right\} \f_r\big(q^{1/m}/\zeta\big) \\
\hphantom{\g_r^{(\zeta)}}{}
= \exp\left\{D\big(P\big(q^{1/m}/\zeta\big)^{Q\p_Q}, Q, q-1\big)\right\}\\
\hphantom{\g_r^{(\zeta)}=}{} \times\exp\left\{\sum_{l>0} \Psi^{lm} \big(D_{lmr}\big(Pq^{lQ\p_Q},Q,q\big)\big)/l^2\big(1-q^m\big)\right\} \f_r^{(\zeta)} .\end{gather*}
The transition from the first to the second line uses a version of the Campbell--Hausdorff formula which rewrites
${\rm e}^{A+B/(1-q)}$ as ${\rm e}^C{\rm e}^{B/(1-q)}$, where by $A$ and $B$ we denote operators without pole at $q=1$. Since the commutator of $[A,B]$ is divisible by $q-1$, the operator~$C$ also comes out without the pole.

Comparing with the formula for $T_{\g_{rm}^{(1)}}\L_X^{\rm fake}$, we find that the factor ${\rm e}^{B/(1-q)}$ coincides with
\[ \Psi^m \left(\exp\left\{ \sum_{k>0} \Psi^k\big(D_{krm}\big(Pq^{kQ\p_Q},Q,q\big)\big)/k^2(1-q)\right\} \right) .\]
Given that $\f_r^{(\zeta)} \in \Delta_{\zeta} \Psi^m\big(T_{\f_{rm}^{(1)}}\L_X^{\rm fake}\big)\otimes \Lambda$, and using the fact that $\Delta_{\zeta}$ does not depend on~$Q$, we conclude, that ${\rm e}^{B/(1-q)}\f_r^{(\zeta)} \in \Delta_{\zeta} \Psi^m\big(T_{\g_{rm}^{(1)}}\L_X^{\rm fake}\big)\otimes \Lambda$.

As we have discussed, the tangent space $T_{\g_{rm}^{(1)}}\L_X^{\rm fake}$ is preserved by (differential and hence) finite difference operators. However, Adams' operation $\Psi^m$ changes the representation by which such operators act. Namely,
$\Psi^m \big(\big(P_iq^{Q_i\p_{Q_i}}\big)^{1/m}\big) = P_i \big(q^{Q_i\p_{Q_i}}\big)^{1/m}$. Also note that
$\zeta^{-Q_i\p_{Q_i}}$, ac\-ting on $Q^d$ by $\zeta^{-d_i}$, acts trivially on $\Psi^m\big(Q^d\big)=Q^{md}$ when $\zeta$ is an $m$th root of unity. Therefore the operator $D\big(P\big(q^{1/m}/\zeta\big)^{Q\p_Q},Q,q-1\big)$ (nicknamed as $C$ in the
Campbell--Hausdorff formula) acts in $\Psi^m\big(\K^{\rm fake}\big)$ as a legitimate pseudo-differential operator in this new representation, and preserves the subspace $\Psi^m\big(T_{\g_{rm}^{(1)}}\L_X^{\rm fake}\big)$. This implies that the same is true for the subspace $\Delta_{\zeta}\Psi^m\big(T_{\g_{rm}^{(1)}}\L_X^{\rm fake}\big) \otimes \Lambda$
in $\K^{\rm fake}$, and shows that $\g_r^{(\zeta)}={\rm e}^{C}{\rm e}^{B/(1-q)}\f_r^{(\zeta)}$ lies
in this subspace. \end{proof}

\begin{corollary}[string flows]\label{Corollary1} The transformations
\[ \f_r \mapsto \exp\left\{ \sum_{k>0} \Psi^k(\tau_{kr})/k\big(1-q^k\big)\right\} \f_r,\qquad r=1,2,\dots, \]
where $\tau_k\in \Lambda_{+}$, $k=1,2,\dots$, preserve $\L_X$.
\end{corollary}

\section{Explicit reconstruction}\label{section4}

We generalize ``explicit reconstruction'' results from \cite{GiE} and Part VIII. Let us begin with the following proposition.

\begin{proposition}[$\D_q$-module structure]\label{Proposition1} If $\f=(\f_1,\dots,\f_r,\dots)$ lies in $\L_X$, then $(D_1\f_1, \dots, D_r\f_r,\allowbreak \dots)$, where $D_1,D_2,\dots$ is any sequence of finite difference operators, also lies in $\L_X$.
\end{proposition}

\begin{proof} This is a simplified version of the previous arguments. We are given that $\f_r^{(\zeta)}$ pass the tests (i), (ii) of Theorem~\ref{Theorem1}, i.e.,
$\f_r^{(1)} \in \L_X^{\rm fake}$, and for primitive $m$th roots of unity $\zeta\neq 1$,
$\f_r^{(\zeta)}\in \Delta_{\zeta} \Psi^m\big(T_{\f_{rm}^{(1)}}\L_X^{\rm fake}\big)\otimes \Lambda$.
Then $\f_r^{(1)}$ lies in (its own) ruling space $(1-q)T_{\f_r^{(1)}}\subset \L_X^{\rm fake}$, which is $D_r$-invariant. Therefore, for $\g_r:=D_r\f_r$, we find that $\g_r^{(1)}$ lies in the same ruling space as $\f_r^{(1)}$, and the tangent spaces $T_{\g_r^{(1)}}$ and $T_{\f_r^{(1)}}$ coincide as well. Furthermore,
\[ \g_r^{(\zeta)} = D_r\big(P\big(q^{1/m}/\zeta\big)^{Q\p_Q}, Q, q^{1/m}/\zeta\big) f_r^{(\zeta)} \in
\Delta_{\zeta} \Psi^m\big(T_{\f_{rm}^{(1)}}\big)\otimes \Lambda ,\]
since the operator here commutes with $\Delta_{\zeta}$ and preserves the subspace $\Psi^m\big(T_{\f_{rm}^{(1)}}\big)$ in $\Psi^m\big(\K^{\rm fake}\big)$ (as was discussed in the proof of Theorem~\ref{Theorem3}). \end{proof}

\begin{remark} Proposition~\ref{Proposition1} together with Corollary~\ref{Corollary1} give another proof of Theorem~\ref{Theorem2}. Namely, it is easy to check that $\t={\mathbf 0}$ is a critical point of $\F_{pt}$, which implies that the dilaton vector $(1-q)(\1,\1,\dots)$ lies in $\L_{pt}$. Then, applying Theorem~\ref{Theorem4} with $D_r=\t_r\big(q,q^{-1}\big)$ (multiplication by a Laurent polynomial),
we conclude that $(1-q)\K^{\infty}_{+}$ lies in $\L_X$. Finally, applying the string flow from Corollary, we obtain the whole of~$\L_{pt}$. \end{remark}

\begin{theorem}\label{Theorem4} Let the degree-$2$ classes $p_1,\dots, p_s$ generate the cohomology algebra $H^{\rm even}(X;\QQ)$, and hence the line bundles $P_1,\dots,P_s$ multiplicatively generate $K^0(X)\otimes \QQ$. Pick Laurent monomials $P^{m_{\a}}$ to form a linear basis in $K^0(X)\otimes \QQ$. Let $c_{\a,r}\in \Lambda \big[q,q^{-1}\big]$, $r=1,2,\dots$, be arbitrary\footnote{For consistency with our earlier understanding of the domain of $\F_X$ we should assume that $c_{\a,r}$ are $\Lambda_{+}$-close to $1$.} scalar Laurent polynomials in $q$, and $\tau_{\a,k} \in \Lambda_{+}$, $k=1,2,\dots$, be arbitrary ``small'' parameters. Suppose that $\f=(\f_1,\dots,\f_r,\dots)$ is a point in $\L_X$, and let $\f_r=\sum_d f_{r,d} Q^d$ be the $Q$-expansion of $\f_r$ $($i.e., the coefficients $f_{r,d}$ do not depend on~$Q)$. With this notation, the following family $\g =(\g_1,\dots, \g_r,\dots)$ of points in $\K^{\infty}$ lies in $\L_X$ and parameterizes the whole of it
 \[ \g_r:= \sum_d f_{r,d} Q^d \exp\left\{\sum_{k>0} \sum_\a \Psi^k(\tau_{\a,rk}) P^{km_a}q^{k(m_\a,d)}/k\big(1-q^k\big)\right\}
 \sum_{\a} c_{\a,r}(q) P^{m_\a}q^{(m_\a,d)}.\]
\end{theorem}

Here $(m, d)=\sum_i m^id_i$ denotes the value of $\sum_i m^i p_i \in H^2(X)$ on
$d\in H_2(X)$.

\begin{proof} Let us assume first that $\tau_{\a,k}$ are free $\lambda$-algebra generators added to the ground ring $\Lambda$. We apply Theorem~\ref{Theorem3} and Proposition~\ref{Proposition1} over such an extended ground ring $\tilde{\Lambda}$ to conclude that the following family $\tilde{\g} = (\dots, \tilde{\g}_r,\dots)$ of points
lies in $\tilde{\L}_X$:
\[ \tilde{\g}_r:=\left(\sum_\a c_{\a,r}\big(Pq^{Q\p_Q}\big)^{m_\a}\right)
\exp\left\{\sum_{k>0} \sum_\a \Psi^k(\tau_{\a,rk}) \big(Pq^{Q\p_Q}\big)^{km_\a}/k\big(1-q^k\big)\right\} \f_r.\]
Since $\tau_{\a,k}$ do not depend on $Q$, the operators in the exponent commute
and can be computed explicitly on the monomials $Q^d$. Namely, $\big(q^{Q\p_Q}\big)^{km_{\a}}Q^d =
q^{k(m_\a,d)}Q^d$. After that the operators on the left can be applied explicitly
too. The resulting expression matches the one for~$\g_r$ in the formulation of Theorem~\ref{Theorem4}. Now we can invoke the change of the ground ring $\tilde{\Lambda}\to \Lambda$ defined as the homomorphism of $\lambda$-algebras induced by the specialization of $\tau_{\a,k}$ to their values in $\Lambda_{+}$. It is easy to see that the descendant potential $\tilde{\F}_X$ turns into~$\F_X$ under this operation, and hence the variety $\tilde{\L}_X$ turns into $\L_X$. Therefore the family $\g$ indeed lies in $\L_X$.

We can check now the same way as we did in the proof of Theorem~\ref{Theorem2} that the projection of this family to $\K_{+}^{\infty}$ along $\K_{-}^{\infty}$ covers the $\Lambda_{+}$-neighborhood of the dilaton vector. For instance, modulo Novikov's variables, taking all $\sum_\a c_{r,\a} P^{m_{\a}} =1$, we have
\[ g_{r,0}=f_{r,0}+\sum_{\a} \tau_{\a,r} P^{m_{\a}} + \text{terms with poles at $q\neq 0,\infty$}, \]
 which gives a family\footnote{We assume that $\L_X\ni \f$ lies in $\K^{\infty}$ in $\Lambda_{+}$-neighborhood of the dilaton vector, and in particular $f_{r,0}=(1-q)\1 \mod \Lambda_{+}$.} $[g_{r,0}]_{+}=(1-q)\1+t$ with any $q$-independent value of $t$. In the other extreme, when all $\tau_{\a,k}=0$,
we have
\[ g_{r,0} = \left(\sum_\a c_{\a,r}(q)|_{Q=0} P^{m_{\a}}\right) f_{r,0},\]
which yields a family $[g_{r,0}]_{+}$ transverse to $(1-q)\1+t$, i.e., in projection along $(1-q)+t$, containing all multiples of $(1-q)$. Our claim follows from this by virtue of the formal implicit function theorem. \end{proof}

\section{Derivation of Theorem~\ref{Theorem1}}\label{section5}

The higher genus formula for the total descendant potential $\D_X$ of permutation-equivariant quantum K-theory of $X$ was obtained in Part IX
by generalizing to the higher genus case the approach from \cite{GiTo},
where a recursion relation for the J-function of quantum K-theory was derived from the virtual Kawasaki-RR formula on moduli spaces of genus~0 stable maps.
Rather than specializing the same method to the genus-$0$ case, we intend
here to extract the genus-$0$ information encoded by the generating function
$\F_X$ from the formula for~$\D_X$. However, we need to warn the reader, that such a derivation may in the end provide few advantages over the first, more direct method.

{\bf A.} Let us first recall from Part~IX those aspects of the adelic formula for $\D_X$ which are relevant for our goal of extracting from it
the genus-0 information. According to the definition of $\D_X$ found in Part IX, the genus 0 descendant potential $\F_X$ enters the total descendant potential $\D_X$ in the form\footnote{This expression already exhibits the main difficulty of our task in comparison with the ``usual'' quantum mechanics. In, say, quantum cohomology theory, the analogous genus expansion would have the form $D_X={\rm e}^{\h^{-1}\F_X+\O(1)}$, so that the quasi-classical limit is simply extracted as $\F_X=\lim\limits_{\h=0} \h \log \D_X$.}
\[ \D_X = \exp\left\{ \sum_{k>0} \h^{-k}\Psi^k(\F_X(\t_k, \t_{2k},\dots,\t_{rk},\dots))/k +
 \O (1)\right\},\]
where $\O(1)$ denotes all terms weighted by non-negative powers of~$\h$, and
containing contributions of higher genus curves.

The adelic expression given for $\D_X$ in Part~IX is built from the {\em adelic tensor product}
\[ \und{\D}_X \big(\big\{ t_r^{(\zeta)} \big\}, \h, Q\big):= \bigotimes_{M=1}^{\infty} \D_{X/\ZZ_M}^{tw} \left( \Big\{ t_{r(\zeta)}^{(\zeta)}/\sqrt{\h^{r(\zeta)}} \Big\}_{\zeta\colon \zeta^M=1} , 1, Q^M\right),\]
where $\D_{X/\ZZ_M}^{tw}\big(\big\{ t^{(\zeta)}_{r(\zeta)} \big\}, \h, Q\big)$ is the total descendant potential of a certain twisted fake quantum K-theory of the orbifold target space~$X/\ZZ_M$. Recall that the inputs $t_{r(\zeta)}^{(\zeta)} \in \K^{\rm fake}_{+}:=K((q-1))$ of~$\D_{X/\ZZ_M}^{tw}$ are labeled by primitive roots of unity $\zeta$ of {\em orders} $m(\zeta)$ dividing $M$, while the subscript $r(\zeta):=M/m(\zeta)$. This adelic product, subject to the suitable dilaton shift, is considered as a~quantum state $\lan \und{\D}_X\ran$ in the quantization of the adelic product of symplectic loop spaces\footnote{The details about the symplectic loop spaces (both global and adelic), and their symplectic forms are found in Part IX, but we quote later some of these details as they become necessary.}
\[ (\und{\K}^{\infty},\und{\Omega}^{\infty}):=\prod_{r=1}^{\infty}\big(\und{\K}^{(r)},
\und{\Omega}^{(r)}\big).\]
Likewise, the dilaton-shifted function $\D_X$ is considered as a quantum state $\lan \D_X\ran$ in the quantization of the symplectic loop space $(\K^{\infty},\Omega^{\infty})$, where $\K^{\infty}$ consists of sequences $\f=(\f_1,\dots,\f_r,\dots)$ of $K$-valued rational
functions of $q$.

The expression for $\lan \D_X \ran$ in terms of
$\lan \und{\D}_X\ran$ is induced by the symplectic {\em adelic map} $\und{\ }\colon (\K^{\infty},\Omega^{\infty})\! \allowbreak \to \big(\und{\K}^{\infty},\und{\Omega}^{\infty}\big)$. To a sequence $\f$ of rational functions rational function of~$q$, it associates a~collection $\big\{ \und{\f}_r^{(\zeta)} \big\}$ of $K$-valued Laurent series in~$(q-1)$, where $\und{\f}_r^{(\zeta)}$ is the Laurent series expansion near $q=1$ of~$\Psi^r\big(\f_r\big(q^{1/m(\zeta)}/\zeta\big)\big)$. Note that $\und{\f}_r^{(\zeta)}$ is related by~$\Psi^r$ to what we denoted earlier in this text by~$\f_r^{(\zeta)}$.

\looseness=-1 The relation between the quantization spaces
induced by the symplectic adelic map involves the change from the {\em standard} Lagrangian polarization $\und{\K}^{\infty}_{\pm}$ (using which the adelic product $\otimes_M \D_{X/\ZZ_M}^{tw}$ is lifted to a function on $\und{\K}^{\infty}$) to the {\em uniform} polarization, in which the negative space coincides with adelic image of $\K^{\infty}_{-} \subset \K^{\infty}$. After the polarization's change, $\lan \D_X\ran$ is identified with the restriction
of $\lan \und{\D}_X\ran$ to the image of $\K^{\infty}_{+}=\K^{\infty}/\K^{\infty}_{-}$ in $\und{\K}^{\infty}_{+}$ under the adelic map.

The whole expression has the structure of Wick's formula of summation over graphs, where the factors in $\otimes_M \lan \D_{X/\ZZ_M}^{tw}\ran$ represent
different types of vertices, and the edge propagators come from the change of polarization. The latter is equivalent to the application of
${\rm e}^{\oplus_{r=1}^{\infty} \h^r \nabla_r/2}$, where~$\nabla_r$ is a certain 2nd order
differential operator acting only through the inputs
$t_r^{(\zeta)}/\sqrt{\h^r}$ of~$\und{\D}_X^{\infty}$ of the same level~$r$.

{\bf B.} Our goal is to extract from this expression a formula for the genus-$0$ part of $\log \D_X$, i.e., the terms with negative powers of~$\h$:
\[ \sum_{k>0} \h^{-k}\Psi^k(\F_0(\t_k,\t_{2k},\dots,\t_{rk},\dots))/k.\]

Recall that each $\D_{X/\ZZ_M}^{tw}$ has the form of the genus expansion $\D (\t) = {\rm e}^{\sum_g \h^g \F_g(\t)}$, where~$\F_g$, which counts contributions of genus-$0$ stable maps to the orbifold $X/\ZZ_M$, is homogeneous of degree $2-2g$ in $\t$ {\em after} the dilaton shift. Contributions to genus-0 part of $\log \D_X$ come, however, from mappings of those connected orbicurves to $X/\ZZ_M$ whose ramified $\ZZ_M$-cover is rational (though not necessarily connected). If $\hat{\gs}_i$, $i=1,\dots, n$, are marked points of ramification index $m_i=M/r_i$ on a connected orbicurve $\hat{\Sigma}$, then Hurwitz' formula for the Euler characteristic of the covering curve $\Sigma$ gives
\[ \eu (\Sigma) = M (2-2g-n)+\sum_{i=1}^n r_i =
M \left( 2-2g - \sum_i \left(1-\frac{1}{m_i}\right)\right).\]
For $\eu (\Sigma)$ to be positive, we must have $g=0$, all but at most $3$ of $m_i$ equal to~$1$, and those~$3$ (call them~$a$, $b$, $c$) satisfy $1/a+1/b+1/c>1$.
Of course, this singles out the $ADE$-orbifold structures on~$\CC P^1$. However, it is easy to check that in the cyclic group $\ZZ_M$ there are no elements of
orders~$a$,~$b$,~$c$ whose product is the identity, unless it is $A_{m-1}$-case: one of $a$, $b$, $c$ equals $1$, and the other~$2$ are equal to~$m$. In other words, if $\Sigma \to \CC P^1$ is a~rational ramified $\ZZ_M$-cover, then base $\CC P^1$ carries~2 ramification points of the same order~$m$, $\Sigma$ consists of $r=M/m$ disjoint copies of~$\CC P^1$, and the generator $h\in \ZZ_M$ acts on~$\Sigma$ by cyclically permuting the $r$ copies in such a way that~$h^r$ acts on each of them as the multiplication by a primitive $m$th root of unity $\zeta$ at one of the ramification points (and hence by~$\zeta^{-1}$ at the other).

We can now extract from the logarithm of the adelic product the vertex contribution into the genus-0 part of $\log \D_X$. It has the form
\[ \sum_{M=1}^{\infty} \left( \F_{X/\ZZ_M}^{tw}\left(\frac{t_M^{(1)}}{\sqrt{\h^M}}\right) + \frac{1}{2}\sum_{\stackrel{\zeta^M=1}{\zeta\neq 1}} \left\lan {\rm d}^2\F^{tw}_{X/\ZZ_M}\left(\frac{t_M^{(1)}}{\sqrt{\h^M}}\right)
\frac{t^{(\zeta)}_{r(\zeta)}}{\sqrt{\h^{r(\zeta)}}}, \frac{t^{(\zeta^{-1})}_{r(\zeta^{-1})}}{\sqrt{\h^{r(\zeta^{-1})}}} \right\ran \right) .\]
In this expression we use a number of {\em ad-hoc} conventions. By $\F_{X/\ZZ_M}^{tw}$ we denote here the genus-$0$ descendant potential of fake twisted quantum K-theory of $X/\ZZ_M$. It is a function of
$\big\{ t_{r(\zeta)}^{(\zeta)} \big\}$, where $\zeta$ is an $M$th root of unity of primitive order $m(\zeta)$, and $r(\zeta)=M/m(\zeta)$. The input~$t_M^{(1)}$
corresponds to the unramified sector ($\zeta=1$). In the left summand, we assume that the inputs from all other sectors are set to~$0$. In the previous discussion, this term represents contributions of totally unramified (and hence trivial) $\ZZ_M$-covers $\Sigma \to \CC P^1$. The right sum represents covers ramified at~$2$ points. Each term comes from the quadratic differential of $\F_{X/\ZZ_M}^{tw}$, and is a~bilinear form in $t_{r(\zeta)}^{(\zeta)}$ and $t_{r(\zeta^{-1})}^{(\zeta^{-1})}$.
The coefficient matrix of the bilinear form consists, of course, of the 2nd derivatives of $\F_{X/\ZZ_M}^{tw}$ in the direction of the appropriate sectors (corresponding to~$\zeta$ and~$\zeta^{-1}$), and evaluated at the input $t_M^{(1)}/\sqrt{\h^M}$ (as indicated in the formula), while {\em all inputs $t_{r(\eta)}^{(\eta)}$ from ramified sectors $(\eta^M=1$, $\eta \neq 1)$ are set to zero}. The linear operator ${\rm d}^2\F_{X/\ZZ_M}^{tw}\big(t_M^{(1)}/\sqrt{\h^M}\big)$ acts from
the space $\big(\K_r^{(\zeta)}\big)_{+}\ni t_r^{(\zeta)}$ to $\big(\K_r^{(\zeta)}\big)_{-}$ which is identified by the symplectic pairing with the dual of $\big(\K_r^{(\zeta^{-1})}\big)_{+}\ni t_r^{(\zeta^{-1})}$.

Since the (dilaton-shifted) functions $\F_{X/\ZZ_M}^{tw}$ have homogeneity degree~$2$, we can rewrite the above sum by powers of~$\h$:
\[ \und{\F}:=\sum_{r=1}^{\infty} \h^{-r} \left( \F_{X/\ZZ_r}^{tw}(t^{(1)}_r) +\frac{1}{2}\sum_{\zeta\neq 1}
\left\lan {\rm d}^2 \F^{tw}_{X/\ZZ_{m(\zeta)r}} \big(t_{m(\zeta)r}^{(1)}\big) t_r^{(\zeta)}, t_r^{(\zeta^{-1})}\right\ran \right) .\]

{\bf C.} The above vertex contribution should be combined with the edge operators $\nabla_r$ (propagators). There are two ways to proceed to our goal of extracting the genus~0 data from the adelic expression. One is to retain in $\log \big({\rm e}^{\oplus_r \h^r \nabla_r/2} {\rm e}^{\und{\F}}\big)$ all terms of negative order in~$\h$. This leads to summation over rooted trees. The other is to introduce the adelic version $\und{\L}^{\infty}_X \subset \big(\und{\K}^{\infty}, \und{\Omega}^{\infty}\big)$ of the variety~$\L_X$ using~$\und{\F}_X$ and the standard polarization $\und{\K}_{\pm}^{\infty}$.

Let us start with the first way. Put $\und{\F}=\sum_r \h^{-r}\und{\F}_r$, i.e.,
denote by $\und{\F}_r$ the sum of the terms weighted by the $r$th power of $\h^{-1}$, and consider each $\und{\F}_r$ as a family of functions of
$t_r =\big\{ t_r^{(\zeta)} \big\}$ depending on $t_{mr}^{(1)}$ with $m>1$ as parameters.
Let us denote the parameters in $\und{\F}_r$ by $t_{m,r}^{(1)}$ (thereby making the parameters with different~$r$ independent of each other). The differential operator $\nabla_r$ acts only through the variables~$t_r$. Therefore
\[ \log \big( {\rm e}^{\h^r\nabla_r/2} {\rm e}^{\h^{-r}\und{\F}_r} \big) = \h^{-r} \F_r + \O\big(\h^0\big),\]
where $\F_r$ is obtained by Wick's summation over trees from vertex contribution $\und{\F}_r$ and edge propagator~$\nabla_r$. It still depends on~$t_{m,r}^{(1)}$ with $m>1$ as parameters.

Let us now employ our second approach, based on the families of Lagrangian varieties in adelic symplectic loop spaces.

We have
\[ \und{\F}_r = \F_{X/\ZZ_r}^{tw}\big(t^{(1)}_r\big) +\frac{1}{2}\sum_{\zeta\neq 1}
\big\lan {\rm d}^2 \F^{tw}_{X/\ZZ_{m(\zeta)r}} \big(t_{m(\zeta),r}^{(1)}\big) t_r^{(\zeta)}, t_r^{(\zeta^{-1})}\big\ran,\]
where $t_{m,r}^{(1)}$ are viewed as parameters. Given their values, the graph of
${\rm d}\und{\F}_r$ is a Lagrangian submanifold $\und{\L}_r$ in the adelic loop space
$\und{\K}^{(r)}$, which is the Cartesian product\footnote{More precisely, it is the subset in this product, consisting of sequences of Laurent series in $q-1$ with all but finitely many terms having no pole when considered modulo any finite term of the filtration in~$\Lambda$.}
\[ \und{\K}^{(r)}:=\prod_{\text{roots of unity $\zeta$}}\ \K_r^{(\zeta)},\]
where each $\K_r^{(\zeta)}$ is a copy of $\K^{\rm fake}:=K((q-1))$. It is equipped with the Cartesian product symplectic form
\[ \und{\Omega}^{\infty}(\f,\g)=\sum_{\zeta} \frac{1}{m(\zeta)} \sum_{r=1}^{\infty} \frac{1}{r} \Res_{q=1}\big(\und{\f}_r^{(\zeta)}\big(q^{-1}\big),\und{\g}_r^{(\zeta^{-1})}(q)\big)^{(r)} \frac{{\rm d}q}{q},\]
where the twisted Poincar\'e pairing $(\cdot ,\cdot )^{(r)}$ on $K:=K^0(X)\otimes \Lambda$ is characterized by
\[ \big(\Psi^r a, \Psi^r b\big)^{(r)} = r \Psi^r(a,b).\]
Note that the symplectic form pairs the Lagrangian subspace $\big(\K_r^{(\zeta)}\big)_{\pm}$ with $\big(\K_r^{(\zeta^{-1})}\big)_{\mp}$ (here $+/-$ refer to ``power series/principal part'' respectively), and the standard polarization is $\und{\K}^{\infty}_{\pm}:=\prod_{\zeta, r} \big(\K_r^{(\zeta)}\big)_{\pm}$.

Based on the explicit form of $\und{\F}_r$, for each value of the parameters, we have $\und{\L}_r = \prod_{\zeta} \L_r^{(\zeta)}$,
where $\L_r^{(1)}\subset \K_r^{(1)}$ is the graph of ${\rm d}\F_{X/\ZZ_r}^{tw}(t_r^{(1)})$, and for $\zeta\neq 1$, each $\L_r^{(\zeta)} \subset \big(\K_r^{(\zeta)}\big)_{+}\oplus \big(\K_r^{(\zeta^{-1})}\big)_{-}$ is the graph of the self-adjoint linear map
$t_r^{(\zeta)} \mapsto {\rm d}^2 \F^{tw}_{X/\ZZ_{m(\zeta)r}} \big(t_{m(\zeta),r}^{(1)}\big) t_r^{(\zeta)}$.

Our reasoning based on summation over trees implies that {\em replacing $\und{\F}_m$ with $\F_m$ can be considered as the change of the standard polarization into the uniform one without changing the Lagrangian submanifold $\und{\L}_r \subset \und{\K}^{(r)}$ which represents the graph of the differential.}

{\bf D.} Using the above conclusion, we now aim to extract terms of weight $\h^{-1}$ from the whole expression. For this we will observe how all $\nabla_{mr}$ (that they carry the weights $\h^{mr}$) intertwine $\F_r$ with $\F_{mr}$ (acting through the parameters $t_{m,r}^{(1)}$ on the former, and through the variables $t_{mr}^{(\zeta)}$ on the latter). We claim that {\em the terms of the total weight $\h^{-r}$ are obtained by Wick's summation over {\em rooted} trees, where the root vertex is represented by~$\F_r$, all other vertices correspond to~$\F_{mr}$ with $m>1$, and the edges, connecting vertices of lower {\em level} $mr$ with higher level~$mlr$, come from the propagators $\nabla_{mlr}$.}

\begin{figure}[htb]\centering
\includegraphics[scale=0.85]{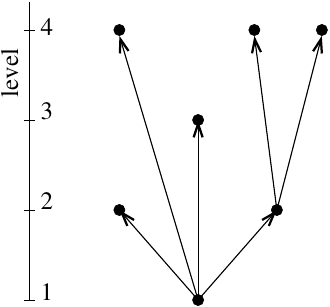}
\caption{A rooted tree.}\label{Fig1}
\end{figure}

This is a purely combinatorial statement. Orient the edges from lower level to higher (Fig.~\ref{Fig1}). For a tree just described, each non-root vertex of level~$mr$ comes with the weight $\h^{-mr}$ and a~unique {\em entering} edge weighted by the inverse factor~$\h^{mr}$. Therefore the whole rooted tree contributes with the weight $\h^{-r}$ of the root vertex. Conversely, suppose that a tree has $k+1>1$ ``roots'' (i.e., vertices with no entering edges). Then there are $k$ edges in excess of those entering non-root vertices and canceling their weights as above. Moreover, since the part of the tree below any fixed level must remain cycle-free, the weights $-n_i$, $n_0\leq \dots \leq n_k$ of the roots are majorated by the weights $m_1\leq \dots \leq m_k$ of the excess edges this way: $n_i<m_i$. In effect, due to the divisibility properties of the weights, $2n_i\leq m_i$. Therefore $n_0+\cdots+n_k \leq m_1+\cdots + m_k$, and hence the total weight of the tree is non-negative. Finally, adding cycle-generating edges to a tree (even with one root) makes the weight of the resulting graph non-negative.

{\bf E.}
Let us examine now the summation over trees with one root of level
$r=1$. It affects the root contribution $\F_1$ by shifting the values of the parameters $t=\big(t_{1,1}^{(1)},\dots,t_{m,1}^{(1)},\dots\big)$. This follows from Taylor's formula
\[ \F_1(t+y)=\sum_{\a=(\a_1,\dots,\a_m,\dots)} \partial^{\a}\F_1(t)\frac{y^\a}{\a!},\]
where $y_m$ signifies the component of the derivation $\nabla_m \F_m$ (recall that we view~$\nabla_m$ as a bi-derivation) in the direction of the variables~$t_{m,1}^{(1)}$ (and dropping the components in the directions of~$t_m^{(\eta)}$ with $\eta\neq 1$). Of course, the variables $t_m^{(\zeta)}$ of $\F_m$ together with partial derivatives of~$\F_m$ in these variables (computed at certain values of the parameters~$t_{l,m}^{(1)}$) form a point on the graph of the differential of~$\F_m$ (at these values of the parameters), which is a Lagrangian submanifold in the symplectic loop space $\und{\K}^{(m)}:=\prod_{\zeta} \K_m^{(\zeta)}$. Let us denote this point by $\und{f}_m$, and the Lagrangian submanifold by~$\und{\L}_m$. Since $\nabla_m$ generates the transition between the standard and uniform polarizations on~$\und{\K}^{(m)}$,
we conclude that {\em the shifted value $t_m^{(1)}+y_m$ of the parameter $t_{m,1}^{(1)}=t_m^{(1)}$ can be described as the component in the untwisted sector of the projection $[\und{f}_m]_{+}$ of $\und{f}_m\in \und{\L}_m\subset \und{\K}^{(m)}$ to $\und{\K}^{(m)}_{+}$ with respect to the {\em uniform} polarization} (while with respect to the standard one, such projection is~$t_m$ by the very definition).

The same is true for the terms weighted by $\h^{-r}$, which come from rooted trees with the root at the level $r$: the effect of all propagators $\nabla_{mlr}$ on $\F_r$ consists in shifting the values of the parameters $t_{rm}^{(1)}$ into $[\und{f}_{rm}]_{+}^{(1)}$ where $\und{f}_{rm}\in \und{\L}_{rm}\subset \und{\K}^{(rm)}$ represents the differential ${\rm d}\F_{rm}$ at the point $\big\{ t_{rm}^{(\eta)} \big\}$ and at the appropriate values of the parameters $t_{rlm}^{(1)}$. The subscript in $[\dots]_{+}^{(1)}$ refers to the projections along the negative space of the {\em uniform} polarization in $\und{\K}^{(rm)}$, while the superscript indicates the component in the untwisted sector $\eta=1$. Note that the ``appropriate'' values of the parameters~$t_{rlm}^{(1)}$ are likewise determined by the differentials of $\F_n$, where~$n$ runs multiples of~$rlm$, so that this description has the form of an infinite recursion relation.

{\bf F.} Our next goal is to identify $\L_r^{(\zeta)}$ in terms of the genus-0 descendant
potential $\F_X^{\rm fake}$ of the fake (non-twisted) quantum K-theory of $X$.
We denote by $\L_X^{\rm fake}\subset \K^{\rm fake}$ the overruled Lag\-rangian cone (see \cite{Co, GiF, GiTo}) representing the graph of the differential of $\F_X^{\rm fake}$ with respect to the Lagrangian polarization $\K^{\rm fake}_{\pm}$ on $\big(\K^{\rm fake}, \Omega^{\rm fake}\big)$ (and dilaton-shifted by $(1-q)\1$). In the symplectic space $\K_r^{(1)}$ (it is equipped with the symplectic form based on the twisted pairing~$(\cdot, \cdot)^{(r)}$), there lies the Lagrangian cone $\Psi^r\big(\L_X^{\rm fake}\big)$. In fact $\Psi^r\big(\L_X^{\rm fake}\big)=\L_r^{(1)}$.

To prove this, let us express in terms of $\F_X^{\rm fake}$ the function $\F_{X/\ZZ_r}^{tw}$ at the unramified input $t^{(1)}_r=t$ and with all inputs from ramified sectors set to zero. We claim that\footnote{$\Psi^r$ is invertible in {\em fake} K-theory.}
\[ \F_{X/\ZZ_r}^{tw}\big(t, Q^r\big) = r^{-1}\Psi^r\big(\F_X^{\rm fake}\big(\Psi^{1/r}t, Q\big)\big).\]
Indeed, let $\M$ temporarily denote a moduli space, $X_{0,n,d}$, of genus-0 stable maps to~$X$, so that its contribution to $\F_X^{\rm fake}(t)$ is given by
\[ \chi^{\rm fake}(\M; V) = \int_{[\M]} \ch (V) \td (T_{\M}),\]
where \looseness=-1 $V=\otimes_{i=1}^n \big(\sum_m\ev_i^*(t_m) L_i^m\big)$, and $T_{\M}$ is the virtual tangent bundle to $\M$. The same moduli space parameterizes
stable maps to $X/\ZZ_r$, which are {\em trivial} $\ZZ_r$-covers of stable maps
to~$X$. Its contribution to $\F_{X/\ZZ_r}$ (non-twisted) equals $r^{-1}\chi^{\rm fake}(\M; V)$, where the factor $r^{-1}$ comes from $\ZZ_r$-symmetries of the covers.
The contribution into $\F_{X/\ZZ_r}^{tw}$ is different, and is equal to
\[ r^{-1}\chi^{\rm fake}\left(\M; V \otimes \frac{\Eu(T_{\M})}{\Eu(\Psi^r(T_{\M}))}\right).\]
The twisting fraction $\Eu(T_{\M})/\Eu\big(\Psi^r(T_{\M})\big)$ reflects the difference between deformations of a~trivial $\ZZ_r$-cover as stable maps to $X/\ZZ_r$ and as a stable map to~$X$ of a disconnected curve with~$r$ components. For each K-theoretic Chern root $L$ of the virtual tangent bundle, we have (putting $x=c_1(L)$):
\[\td (L) \left(\ch \frac{\Eu (L)}{\Eu (L^r)}\right) = \frac{x}{1-{\rm e}^{-x}} \frac{1-{\rm e}^{-x}}{1-{\rm e}^{-rx}} = \frac{1}{r} \td \big(L^r\big).\]
Therefore in cohomological terms, the contribution into $\F_{X/\ZZ_r}^{tw}$
reads
\[ r^{-1-\dim_{\CC} T_{\M}}\int_{[\M]} \ch(V) \td \big(\Psi^r(T_{\M})\big). \]
 The integrand here can be rewritten as $\Psi^r \big(\ch \big(\Psi^{1/r}(V)\big) \td (T_{\M})\big)$. Note that $\Psi^r$ acts on cohomology as multiplication by $r^{degree/2}$. Since integration over the virtual fundamental cycle~$[\M]$
 picks only the terms of degree $2\dim_{\CC} T_{\M}$, the result comes out the same as
 \[ r^{-1} \int_{[\M]} \ch \big(\Psi^{1/r}(V)\big) \td (T_{\M}) = r^{-1} \chi^{\rm fake}\big(\M ;
 \Psi^{1/r}(V)\big).\]
 Taking into account that the polarizations in $\K_1^{(1)}=\K^{\rm fake}$ and $\K_r^{(1)}$ are also related by the operation $\Psi^r$, we conclude\footnote{The coefficient $r^{-1}$ is absorbed by the fact the symplectic form in $\und{\K}^{(r)}$ satisfies
$\und{\Omega}^{(r)} \big(\Psi^r(f),\Psi^r(g)\big) = r^{-1}\Psi^r \big(\Omega^{\rm fake}(f, g)\big)$.
Indeed, due to degree-2 homogeneity of $\F_X^{\rm fake}$, for a point $f\in \L_X^{\rm fake}$, we have $\F_X^{\rm fake}([f]_+) = \Omega^{\rm fake}(f,[f]_{+})/2$, and hence
$\und{\Omega}^{(r)}\big(\Psi^r(f),\Psi^r([f]_{+})\big)/2=r^{-1}\Psi^r
\big(\F_X^{\rm fake}([f]_{+})\big)$.} that
 the graph $\L_r^{(1)}\subset \K_r^{(1)}$ of the differential of $\F_{X/\ZZ_r}^{tw}$
 in the unramified sector (i.e., after setting $t_r^{(\zeta)}=0$ for all $\zeta\neq 1$) coincides with $\Psi^r\big(\L_X^{\rm fake}\big)$.

{\bf G.} Having identified $\L_r^{(1)}$, let's proceed to $\L_r^{(\zeta)}$ with $\zeta\neq 1$. From Section~3 of Part~IX, we have (in dilaton-shifted notation)
 \[ \F_{X/\ZZ_M}\big(\big\{ \t^{(h)} \big\}_{h\in \ZZ_r}\big) = \sum_{\chi \in \operatorname{Repr} (\ZZ_r)}\F_X^{\rm fake}\left(\frac{1}{M} \sum_{\chi} t^{(h)} \chi(h) \right) .\]
 The quadratic differential of $\F_{X/\ZZ_M}$ at an unramified input is
 \[ \frac{1}{2M} \sum_{h\neq 1} \big\lan {\rm d}^2\F_X^{\rm fake}\big(t^{(1)}\big) t^{(h)}, t^{(h^{-1})}\big\ran ,\]
 where we took into account that dilaton-shifted 2nd derivatives of $\F_X^{\rm fake}$ are homogeneous of degree~$0$. The Lagrangian subspace generated by the quadratic form $\frac{1}{2} \big\lan {\rm d}^2\F_X^{\rm fake}(t) u, u\big\ran$, i.e., the graph of the linear map $u\mapsto {\rm d}^2\F_X^{\rm fake}(t) u$, is the tangent space~$T$ to $\L_X^{\rm fake}$ at the point corresponding to ${\rm d}\F_X^{\rm fake}$ at the input~$t$. Now the transformation from (non-twisted) $\F_{X/\ZZ_M}$ to $\F_{X/\ZZ_M}^{tw}$ described in Section~4 of Part~IX shows that the graph $\L_{r(\zeta)}^{(\zeta)}$ of the linear map $u^{(\zeta)} \mapsto {\rm d}^2\F_{X/\ZZ_M}^{tw}(t) u^{(\zeta)}$ (defined by 2nd derivatives of $\F_{X/\ZZ_M}^{tw}$ in sectors labeled by $\zeta^{\pm 1}$) is obtained from~$T$ by the multiplication operator (where $r=r(\zeta)$, $m=m(\zeta)$)
 \[ \square_{\zeta, r} = \exp\left\{\sum_{k>0}\left(\frac{\Psi^{kr}(T^*_X-1)}{k(1-\zeta^{-k}q^{kr/m})}-
 \frac{\Psi^k(T^*_X-1)}{k(1-q^k)}\right)\right\}.\]
Note that our previous result $\F_{X/\ZZ_M}^{tw}(t)=M^{-1}\F_X^{\rm fake}\big(\Psi^{1/M}(t)\big)$
implies that $\square_{1,M} T = \Psi^M(T)$. (Here $\square_{1,M}$ is $\square_{\zeta,r}$ with $\zeta=1$, $m(\zeta)=1$ and hence $r(\zeta)=M$.) Therefore
\[ \L_{r}^{(\zeta)} = \square_{\zeta, r(\zeta)} \square_{1,M}^{-1} \Psi^M(T)= \Psi^r(\Delta_{\zeta}) \Psi^M(T) = \Psi^r \big( \Delta_{\zeta} \Psi^m(T)\big),\]
where $m=m(\zeta)$, $r=r(\zeta)$, $M=mr$, and
 \[ \Delta_{\zeta}:= \exp\left\{ \sum_{k>0}\left(\frac{\Psi^{k}(T^*_X-1)}{k(1-\zeta^{-k}q^{k/m})} -\frac{\Psi^{km}(T^*_X-1)}{k(1-q^{km})}\right)\right\}.
 \]
 To be more accurate, one needs to tensor the result with $\Lambda$ over $\Psi^M(\Lambda)$ since $\Psi^M$ may not be invertible on the coefficient ring.

{\bf H.} Finally, let us combine all our previous observations on the adelic data with the description of the adelic map which, according to the main theorem of Part IX, induces $\lan\D_X\ran$ from $\lan \und{\D}_X\ran$.
 The term $\F_X$ in $\log \D_X$ is induced this way from terms of weight $\h^{-1}$ in $\log \und{D}_X$ {\em after} the change of standard polarization to the uniform one.

Let $\f =(\f_1,\f_2,\dots,\f_r,\dots)$ be a point on the graph $\L_X\subset \K^{\infty}$ of the differential of $\F_X$. Then for each $r$,
$\und{\f}_r =\big\{ \und{\f}_r\big(q^{1/m(\zeta)}/\zeta\big) \big\}$ lies in $\und{\L}_r=\prod_{\zeta} \L_r^{(\zeta)}$, which is the graph of the differential of $\F_r$ (in fact of one of the functions of the family~$\F_r$, taken at an appropriate value of the parameter) with respect to the uniform polarization.

More concretely, decompose each vector-valued rational function $\f_r$ into the sum of a Laurent polynomial $\t_r\in \K^{(r)}_{+}$ and simple fractions $\pp_r^{(\eta)}$ with poles (or any order) at a root of unity $q=\eta$:
$\f_r = \t_r + \sum_{\eta} \pp_r^{(\eta)}$. Computing the adelic map,
$\und{\f}_r = \big\{ \und{\f}_r^{(\zeta)} \big\}$, for a primitive $m$th root of unity $\zeta$, we find\footnote{Expansions of rational functions into Laurent series near $q=1$ are tacitly assumed on the right.}
\begin{gather*} \und{\f}_r^{(\zeta)} = \Psi^r\big(\pp_r^{(\zeta)}\big(q^{1/m}/\zeta\big)\big) \quad \big(\text{$\K^{(\zeta)}_{-}$-part in the standard polarization}\big) \\
\hphantom{\und{\f}_r^{(\zeta)} =}{} +\Psi^r\left( \t_r\big(q^{1/m}/\zeta\big) + \sum_{\eta\neq\zeta} \pp_r^{(\eta)}\big(q^{1/m}/\zeta\big)\right) \quad \big(\text{$\K^{(\zeta)}_{+}$-part}\big).
\end{gather*}
In particular, if we put $\t_r=0$, these expansions would describe vectors of the negative space in the {\em uniform} polarization. For a general value of~$\t_r$, projections along such vectors yields $\Psi^r\big(\t_r\big(q^{1/m}/\zeta\big)\big)$. This means that $\und{\f}_r$ represents with respect to the {\em uniform} polarization the differential~${\rm d}\F_r$, computed at the input $\big\{ t_r^{(\zeta)}\big\}$ with $t_r^{(\zeta)} =\Psi^r\big(\t_r\big(q^{1/m}/\zeta\big)\big)$. Equivalently: $\und{\f}_r$~represents with respect to the {\em standard} polarization the differential of ${\rm d}\und{\F}_r$, computed at the point where $t_r^{(\zeta)} =\Psi^r\big(\t_r\big(q^{1/m}/\zeta\big) + \sum_{\eta\neq\zeta} \pp_r^{(\eta)}\big(q^{1/m}/\zeta\big)\big)$.
Since the graph of ${\rm d}\und{\F}_r$ is the product $\und{\L}_r=\prod_{\zeta} \L_r^{(\zeta)}$, we conclude that $\und{\f}_r^{(\zeta)} \in \L_r^{(\zeta)}$.

Now recall that $\L_r^{(1)}=\Psi^r\big(\L_X^{\rm fake}\big)$ to conclude that $\f_r^{(1)}=\Psi^{1/r}\big(\und{\f}_r^{(1)}\big)$ lies in $\L_X^{\rm fake}$. This confirms the criterion~(i) in Theorem~\ref{Theorem1}.

Furthermore, to verify the criterion (ii), take $\zeta\neq 1$, and recall that
$\L_r^{(\zeta)}=\Psi^r\big(\Delta_{\zeta}\Psi^m(T)\big)\otimes \Lambda$, where $T$ is some tangent space to $\L_X^{\rm fake}$. We conclude that $\f_r^{(\zeta)}=\Psi^{1/r}\und{\f}_r^{(\zeta)}$ lies
in $\Delta_{\zeta}\Psi^m(T)\otimes \Lambda$.

How is the space $T$ determined? It is $T_f\L_X^{\rm fake}$, where
$\big[\Psi^M(f)\big]_{+}$ ($M=mr$) is the value of the parameter $t_{m,r}^{(1)}$ in the family of~functions $\F_r$. How is that value determined? It is the untwisted component $[\und{\f}_{mr}]_{+}^{(1)}$ of the point on the graph of ${\rm d}\F_{mr}$. As we already know, $\und{\f}_{mr}^{(1)}$ lies in $\Psi^{mr}(\L_X^{\rm fake})$, and hence
$\f_{mr}^{(1)}=\Psi^{1/mr}\big(\und{\f}_{mr}^{(1)}\big)$ lies in $\L_X^{\rm fake}$. This shows
$f=\f_{mr}^{(1)}$ and therefore $T = T_{\f_{mr}^{(1)}}\L_X^{\rm fake}$ as required.

Ultimately we conclude that if $\f=(\f_1,\dots, \f_r,\dots)$ is a point in $\L_X$, then the criteria~(i) and~(ii) of Theorem~\ref{Theorem1} must be satisfied:
 \[ (i)\ \ \f_r^{(1)}\in \L_X^{\rm fake}, \qquad \text{and}\qquad (ii)\ \ \f_r^{(\zeta)}\in \Delta_{\zeta}\Psi^m\big(T_{\f_{rm}^{(1)}}\L_X^{\rm fake}\big)\otimes_{\Psi^m(\Lambda)}\Lambda\]
 for every $r=1,2,\dots$, and every primitive $m$th root of unity $\zeta\neq 1$.

From the computational point of view, this is an infinite system of recursion relations. However, as it is easy to see, modulo any fixed power of the ideal~$\Lambda_{+}$ (or, more explicitly, under a~limited degree and/or the number of marked points on the curves) this yields a finite system of recursion relations which uniquely determines $\f \in \L_X$ from the projection~$\t$ of $\f$ to $\K_{+}^{\infty}$.

Thus, the criteria (i) and (ii) completely characterize points in~$\L_X$ which in a $\Lambda_{+}$-neighbor\-hood of the dilaton vector.

\vspace{-2mm}

\subsection*{Acknowledgements}
\vspace{-1mm}

This material is based upon work supported by the National
Science Foundation under Grant DMS-1611839, by the IBS Center for Geometry
and Physics, POSTECH, Korea, and by IHES, France.

\vspace{-2mm}

\pdfbookmark[1]{References}{ref}
\LastPageEnding

\end{document}